\documentclass[10pt,onecolumn,draftcls]{IEEEtran}

\makeatother
\IEEEoverridecommandlockouts

\usepackage{lipsum}
\usepackage{cite}
\usepackage{amsmath,amssymb,amsfonts}
\usepackage{algorithm,algorithmic}
\usepackage{hyperref}
\hypersetup{hidelinks=true}
\usepackage{textcomp}
\makeatother

\usepackage{graphicx}
\usepackage{amsthm}
\usepackage{tabu}
\usepackage{breqn}
\usepackage{verbatim}
\usepackage{caption}
\usepackage{esvect}
\usepackage{footnote}
\usepackage{siunitx}
\usepackage{color}
\usepackage{url}
\usepackage[makeroom]{cancel}
\usepackage{todonotes}
\newtheorem{theorem}{Theorem}[section]

\newtheorem{proposition}[theorem]{Proposition}
\newtheorem{remark}{Remark}
\newtheorem*{problem statement}{Problem Statement}

\newtheorem{assumption}{Assumption}

\newcommand{\real}{\mathbb{R}}

\newcommand{\lie}{\mathcal{L}}

\newcommand{\Nc}{{\mathcal{N}}}

\newcommand{\Xc}{{\mathcal{X}}}

\DeclareMathOperator*{\argmin}{arg\,min}

\newcommand{\longthmtitle}[1]{\mbox{}\emph{(#1):}}
\newcommand{\setdef}[2]{\{#1 : #2\}}

\newcommand{\norm}[1]{\left\lVert#1\right\rVert}

\begin{document}
\title{Feedback Optimization with State Constraints through Control Barrier Functions}
\author{Giannis Delimpaltadakis$^*$ \quad Pol Mestres$^*$ \quad Jorge Cortés \quad W.P.M.H. Heemels
\thanks{$^*$Equal contribution.}
\thanks{Giannis Delimpaltadakis and W.P.M.H. Heemels are with the Control Systems Technology group, Mechanical Engineering, Eindhoven University of Technology. Pol Mestres and Jorge Cortés are with  the Department of Mechanical and Aerospace
Engineering, University of California San Diego, USA. Emails: \texttt{i.delimpaltadakis@tue.nl, pomestre@ucsd.edu, cortes@ucsd.edu, w.p.m.h.heemels@tue.nl}. \newline \indent This research is partially funded by the European Research Council (ERC) under the Advanced ERC grant PROACTHIS, no. 101055384 and by AFOSR Award FA9550-23-1-0740.}}
\date{}
\maketitle

%
%

\begin{abstract}
    Recently, there has been a surge of research on a class of methods called \emph{feedback optimization}. These are methods to steer the state of a control system to an equilibrium that arises as the solution of an optimization problem. Despite the growing literature on the topic, the important problem of enforcing state constraints at all times remains unaddressed. In this work, we present the first feedback-optimization method that enforces state constraints. The method combines a class of dynamics called \emph{safe gradient flows} with \emph{high-order control barrier functions}. We provide a number of results on our proposed controller, including well-posedness guarantees, anytime constraint-satisfaction guarantees, equivalence between the closed-loop's equilibria and the optimization problem's critical points, and local asymptotic stability of optima.
\end{abstract}

\section{Introduction}
\emph{Feedback optimization} refers to an emerging class of methods addressing the problem of steering the state of a control system to an equilibrium which arises as the solution of an optimization problem; see \cite{hauswirth2021optimization} for a survey. 
%
%
These methods design dynamic feedback controllers, which resemble optimization algorithms, namely \emph{gradient flows}, that drive the system to the desired equilibrium, without explicitly solving the optimization problem. Compared to the feed-forward approach of solving the optimization problem offline and driving the plant to the computed setpoint, feedback optimization enjoys the robustness that comes with feedback and is able to address potentially unknown and time-varying plant dynamics and objective functions. As such, feedback optimization has found a wide range of applications on, e.g., power systems \cite{dall2016optimal, chen2020distributed}, traffic control \cite{GB-JC-JIP-EDA:22-tcns}, smart buildings \cite{belgioioso2021sampled}, communication networks\cite{wang2011control}, etc. Overall, in recent years, research on feedback optimization has surged \cite{dall2016optimal, chen2020distributed, GB-JC-JIP-EDA:22-tcns, KH-JPH-KU:14, belgioioso2021sampled, wang2011control, jokic2009constrained, hauswirth2020timescale, MC-ED-AB:20, brunner2012feedback, LSPL-JWSP-EM:24, YC-LC-JC-EDA:23-csl, haberle2020non, colot2024optimal}. However, although designs with input constraints are available in the literature, a fundamental problem remains unsolved: \emph{enforcing state constraints at all times}. 

%
%
\subsection*{Contributions}
In this work, \emph{we develop a feedback-optimization method for dynamic plants that enforces state and input constraints at all times}; to our knowledge, the first one to achieve so. The method combines so-called \emph{safe gradient flows} (SGFs)~\cite{allibhoy2023control} with \emph{high-order Control Barrier Functions} (CBFs)~\cite{xiao2021high}. 
%
%
The controller dynamics arise as the solution of a quadratic program (QP), which can be solved online. We provide the following results on the proposed dynamics:
\begin{itemize}
    \item We introduce mild conditions that guarantee feasibility and Lipschitz continuity of the controller dynamics' QP, which together guarantee the closed-loop dynamics' well-posedness;
    \item We show that the proposed controller can handle constraints of arbitrary relative degree and renders a subset of the constraint set forward invariant. We also prove that this forward-invariant set contains all critical points\footnote{Points that satisfy the \emph{Karush-Kuhn-Tucker} (KKT) conditions (see \cite{NA-AE-MP:20}), including optima.} of the feedback-optimization problem;
    \item We provide conditions guaranteeing equivalence between the feedback-optimization problem's critical points and the closed-loop's equilibria. Among others, equivalence holds if the considered point lies in the interior of the state constraint set;
    \item We prove local asymptotic stability of local optima, assuming that the optimum lies in the interior of the state constraint set;
    \item We introduce a way to regularize the original optimization problem such that all optimizers are in the interior, thus enabling the aforementioned results on equivalence and stability, while guaranteeing arbitrarily small suboptimality bounds.
\end{itemize}

\subsection*{Related work}
The existing literature on feedback optimization either disregards state constraints, enforces them only asymptotically, e.g. \cite{GB-JC-JIP-EDA:22-tcns, MC-ED-AB:20, brunner2012feedback, LSPL-JWSP-EM:24, YC-LC-JC-EDA:23-csl}, or considers only \emph{static plants}, e.g. \cite{haberle2020non,colot2024optimal}. In the latter two cases, which are equivalent, state constraints can be enforced through input constraints, as there is a one-to-one correspondence between state and input. Applying these methods to dynamic plants may lead to constraint violations in the transients. In contrast, the method proposed here is the first one to enforce state constraints at all times for dynamic plants. Nonetheless, it comes with local stability guarantees and assumes full knowledge of the plant dynamics, compared to other works (e.g., \cite{hauswirth2020timescale}) that guarantee global convergence to critical points or optimizers, and employ minimal information on system dynamics. Whether our method extends to model-free scenarios and can ensure global convergence remains a topic of future research.

SGFs were first introduced in \cite{allibhoy2023control}, as dynamical systems that solve constrained optimization problems. These are modified gradient flows, that employ CBFs, to enforce constraint satisfaction. Further, in \cite{delimpaltadakis2023relationship,delimpaltadakis2024continuous}, it was shown how SGFs are continuous approximations of projected gradient flows. SGFs have been used in
reinforcement-learning to ensure safety in \cite{JF-WC-JC-YS:23-csl,PM-AM-JC:25-l4dc}. In the context of feedback optimization, they have been used as controllers interconnected with dynamic plants in
\cite{colot2024optimal, YC-LC-JC-EDA:23-csl}, where however state constraints are only asymptotically enforced. Our work is the first to combine SGFs with high-order CBFs. In fact, high-order CBFs are the key to enforce state constraints at all times, since these are of relative degree greater than one with respect to the controller dynamics and cannot be handled by standard SGFs interconnected to a plant, as SGFs (and projected gradient flows) manipulate only the controller's vector-field. High-order CBFs enable modifications of the controller's vector field that enforce state constraints.

Finally, while we focus on feedback optimization, there are other methods targeting similar problems, such as extremum seeking (ES; e.g. \cite{chen2025continuous,williams2024semiglobal}) and online convex optimization (OCO; e.g. \cite{li2021online,zhou2023safe, karapetyan2023online, nonhoff2022online}). Like existing feedback-optimization methods, ES methods that address state constraints either cast them as input constraints, implicitly assuming that the plant operates at the steady-state, or consider static plants (which are equivalent approaches); see \cite{chen2025continuous, williams2024semiglobal} and references therein. For differences between feedback optimization and ES the reader is referred to \cite{hauswirth2021optimization}. With respect to OCO, the main differences are that OCO  generally considers discrete-time plants and focuses on finite-horizon regret guarantees, rather than asymptotic convergence to the optimizer.

\section{Preliminaries}\label{sec:prelims}
For $r\in\mathbb{Z}_{>0}$, we denote $[r]:=\{ 1, \hdots, r \}$.  We adopt the following conventions on dimensions: given $g:\real^n\times\real^m\to\real^k$, $(x,u)\mapsto g(x,u)$ continuously differentiable, then $\dfrac{\partial g}{\partial x}(x,u)\in\real^{k\times n}$; given $g:\real^n\to\real$, $x\mapsto g(x)$ continuously differentiable, then $\nabla g(x) = (\frac{\partial g}{\partial x}(x))^\top\in\real^n$. Consider the parametric optimization problem
\begin{equation}\label{eq:parametric_nlp}
    \min_x \Phi(x,\theta) \ \text{s.t.} \ a_i(x,\theta)\geq 0, \ i\in[p],
\end{equation}
where $\Phi,a_i:\real^n\times\real^p\to\real$ and $\theta$ is a parameter. Denote $J(x,\theta):=\{i\in[p]: a_i(x,\theta)=0\}$. The following conditions are commonly used in parametric optimization for \emph{sensitivity analysis} \cite{JL:95}, although often more relaxed conditions suffice:
\begin{itemize}
    \item \emph{Mangasarian-Fromovitz Constraint Qualification} (MFCQ): Problem \eqref{eq:parametric_nlp} satisfies the MFCQ at $(x_*,\theta_*)$ if there exists $q\in\real^n$ such that $\frac{\partial a_i}{\partial x}(x_*,\theta_*)q>0$, for all $i\in J(x_*,\theta_*)$.
    \item \emph{Constant Rank Constraint Qualification} (CRCQ): Problem \eqref{eq:parametric_nlp} satisfies the CRCQ at $(x_*,\theta_*)$ if any subset of $\{\frac{\partial a_i}{\partial x}(x,\theta)\}_{i\in J(x,\theta)}$ remains of constant rank in a neighbourhood of $(x_*,\theta_*)$.
\end{itemize}
If the CRCQ or MFCQ holds at $(x_*,\theta_*)$, the KKT conditions \cite{NA-AE-MP:20} are necessary for optimality of $x_*$, for $\theta=\theta_*$.

\section{Problem Statement}
Consider the control system
\begin{align}\label{eq:system}
    \dot{\xi}(t) = f(\xi(t),\upsilon(t)),
\end{align}
where $f:\real^n\times\real^m\to\real^n$ is the vector field, and $\xi(t)$, $\upsilon(t)$ are the state and control input at time $t$, respectively.

\begin{assumption}\longthmtitle{Globally exponentially stable equilibrium}\label{assum:stability}
    For any constant input signal 
    $\upsilon(t) = u$,
    system \eqref{eq:system} admits a unique equilibrium, which is globally exponentially stable. 
\end{assumption}

This assumption is standard in the literature of feedback optimization and is satisfied in many relevant applications (see, e.g., \cite{hauswirth2021optimization}). It guarantees that controller dynamics can be made slow enough, compared to the plant dynamics, so that stability is not compromised when the controller is interconnected to the plant.  For convenience, we let  $w:\real^m\to\real^n$ denote the map that assigns each input $u\in\real^m$ to the unique equilibrium $w(u)$, i.e.,  $f(w(u),u) = 0$.

We wish to design a controller $\upsilon$
that drives system \eqref{eq:system} to a steady-state $(x_*,u_*)$, which solves an optimization problem:
\begin{equation}\label{eq:steady-state_opti_problem}
    (x_*,u_*)\in\left\{\begin{aligned}\argmin_{x,u} \text{ }&\Phi(x)\\
            \mathrm{s.t.} \text{ } & x=w(u), \, h(x)\geq 0, \, b(u)\geq 0,\end{aligned}
            \right.
\end{equation}
where $\Phi:\real^n\to\real$, $h:\real^n\to\real$ and $b:\real^m\to\real$ are sufficiently smooth functions. The equality constraint $x=w(u)$ is precisely the steady-state constraint. Further, we require the input and state inequality constraints to be respected at all times, that is, for all $t$:
\begin{equation}\label{eq:transient_constraints}
    \begin{aligned}
        (\xi(t),\upsilon(t)) \in \mathcal{X}:=\{(x,u)\in\real^n\times\real^m: \ h(x)\geq 0, \ b(u) \geq 0\}.
    \end{aligned}    
\end{equation}
The above problem, \emph{without the transient state constraints $h(\xi(t))\geq 0$}, has been solved by designing a dynamic controller, with dynamics given by, e.g., a projected gradient flow \cite{hauswirth2021optimization, hauswirth2020timescale} or an SGF \cite{YC-LC-JC-EDA:23-csl}. Nonetheless, incorporating state constraints remains an open problem. This is the problem we aim to solve here.

\section{High-order SGFs for feedback optimization with state constraints}
In the following, we denote $\lie_f h(x,u) = \dfrac{\partial h}{\partial x}(x) \cdot f(x,u)$ and, abusively, $\lie^i_f h(x,u) = \dfrac{\partial \lie^{i-1}_fh}{\partial x} \cdot f(x,u)$, for $i\geq 2$.
We further define $h_i:\real^n\times\real^m\to\real$ for $i\in\mathbb{Z}_{\geq0}$ as:
    \begin{align}\label{eq:h_i}
        h_0(x,u) &:= h(x), \\
        h_i(x,u) &:= \dfrac{\partial h_{i-1}}{\partial x}(x,u) \cdot f(x,u) + \beta h_{i-1}(x,u), \quad i \geq 1, \notag
    \end{align}
where $\beta>0$ is a design parameter. In what follows, we assume $\beta>0$ is fixed. We also define the sets $\mathcal{X}_i:=\{(x,u)\in\real^n\times\real^m: h(x) \geq 0, b(u)\geq 0, h_i(x,u)\geq 0\}$ for $i\in\mathbb{Z}_{\geq 0}$. Note that $\mathcal{X}_i\subseteq \mathcal{X}$, for all $i \in\mathbb{Z}_{\ge 0}$, and $\Xc_0=\Xc$.

\begin{remark}\longthmtitle{Effect of $\beta$ on $\mathcal{X}_i$}\label{rem:tuning beta}
    By increasing $\beta$, the sets $\mathcal{X}_i$ become closer to $\mathcal{X}$.
    Informally, when $\mathcal{X}$ is compact, and $h_i$ and $f$ are continuous, for $\beta$ sufficiently large, we have that $h_1(x,u) \approx \beta h(x)$ and thus $\mathcal{X}_1 \approx \mathcal{X}$. Inductively, for any $i$, we obtain $\mathcal{X}_i \approx \mathcal{X}$, for $\beta$ sufficiently large.
\end{remark}

Towards solving the feedback optimization problem with state constraints, we make the following assumptions:
\begin{assumption}\longthmtitle{Relative degree and regularity}\label{assum:differentiability, regularity and relative degree}
    \begin{enumerate}
        \item There exists $r\in\mathbb{Z}_{>0}$ such that:
        \begin{itemize}
            \item the functions $f$ and $h$ are $(r+1)$-times and $(r+2)$-times continuously differentiable, respectively;
            \item $\mathcal{X}_r \neq \emptyset$; $\frac{\partial \lie^{r}_{f}h}{\partial u}(x,u)\neq 0$, for all $(x,u)\in\mathcal{X}_r$ such that $h_r(x,u)=0$; and $\frac{\partial \lie^{i}_{f}h}{\partial u}(x,u)= 0$, for all $(x,u)\in\mathcal{X}_r$ and all $i<r$.
        \end{itemize} \label{assum_item:relative degree}
        \item The function $b$ is twice continuously differentiable, and $b(u) = 0 \implies \nabla b(u) \neq 0$, for all $(x,u)\in \mathcal{X}$.
        \item $\Phi$ and $w$ are twice continuously differentiable.
    \end{enumerate}
\end{assumption}
The key here is item \ref{assum_item:relative degree}. It requires that the control input appears, after differentiating $h$ along the system $(r+1)$-times, which is basically a relative-degree type of assumption. Note that, under item \ref{assum_item:relative degree}, all functions $\{h_i\}_{i<r}$ depend solely on~$x$, and the function $h_r$ is an $r$-order CBF for system \eqref{eq:system}, cf.~\cite{xiao2021high}. 

We design the controller $\dot{\upsilon}(t) = g_{\epsilon,\alpha,\gamma}(\xi(t),\upsilon(t))$, where:
\begin{equation}\label{eq:sgf_QP}
\begin{aligned}
    &g_{\epsilon,\alpha,\gamma}(x,u)=\left\{\begin{aligned}\argmin_{q}\text{ }& \frac{1}{2}\|q +\epsilon\frac{\partial w}{\partial u}(u)^\top \nabla\Phi(x)\|^2 \\
            \mathrm{s.t.:} \text{ } &\nabla^\top b(u)\cdot q +\alpha b(u)\geq 0\\
            &\begin{aligned}\frac{\partial h_r}{\partial x}(x,u)\cdot f(x,u)+ \frac{\partial h_r}{\partial u}(x,u)\cdot q+\gamma h_{r}(x,u)\geq 0\end{aligned}
            \end{aligned}\right.
\end{aligned}
\end{equation}
and $\epsilon>0$, $\alpha>0$ and $\gamma>0$ are parameters to be designed. The closed-loop dynamics become:
\begin{subequations}
    \begin{align}
        \dot{\xi}(t) &= f(\xi(t),\upsilon(t)),\label{eq:plant-interconnection} \\
        \dot{\upsilon}(t) &= g_{\epsilon,\alpha,\gamma}(\xi(t),\upsilon(t)).
    \end{align}
    \label{eq:interconnected-system}
\end{subequations}
Let us provide intuition behind the design of~\eqref{eq:sgf_QP}. It is known, see e.g.,\cite{hauswirth2021optimization}, that the standard gradient flow $\dot{\upsilon}(t) = -\epsilon\frac{\partial w}{\partial u}(\upsilon)^\top \nabla\Phi(\xi)$, drives the system \eqref{eq:system} to the critical points of the unconstrained version of \eqref{eq:steady-state_opti_problem}, for sufficiently small~$\epsilon$. Inspired by the SGF design proposed in~\cite{allibhoy2023control}, \eqref{eq:sgf_QP} modifies the standard gradient-flow point-wise minimally ($\argmin$), so that the conditions of the QP \eqref{eq:sgf_QP} hold. In the coming section, through CBF arguments, we show that these conditions enforce $b(\upsilon(t))\geq 0$ and $h_r(\xi(t))\geq 0$ for all $t$, assuming $b(\upsilon(0))\geq 0$ and $h_r(\xi(0))\geq 0$, guaranteeing forward invariance of a subset of $\mathcal{X}$. Importantly, employing the high-order CBF $h_r$ is the key to enforce state constraints $h(\xi(t))\geq 0$ by modifying controller dynamics.

\begin{remark}\longthmtitle{Connections with the literature}
    The work~\cite{YC-LC-JC-EDA:23-csl} also proposes a controller based on SGFs \cite{allibhoy2023control}, but only guarantees anytime satisfaction of the input constraints. To address anytime satisfaction of state constraints, the design~\eqref{eq:sgf_QP} includes the constraint on the high-order CBF $h_r$, which significantly complicates the technical analysis, as we illustrate below. Incidentally, results in \cite{allibhoy2023control,delimpaltadakis2023relationship,delimpaltadakis2024continuous} show how the controller in \cite{YC-LC-JC-EDA:23-csl} approximates projected gradient flows, which have also been proposed as a solution to feedback optimization with input constraints (see \cite{hauswirth2020timescale,hauswirth2021optimization}). 
\end{remark}

In what follows, we provide a number of results concerning the proposed dynamics~\eqref{eq:interconnected-system}: feasibility of the QP in \eqref{eq:sgf_QP}, Lipschitz continuity of \eqref{eq:sgf_QP}, constraint satisfaction, relationship between the equilibria of \eqref{eq:interconnected-system} and critical points of \eqref{eq:steady-state_opti_problem}, and stability of minimizers of \eqref{eq:steady-state_opti_problem}. All proofs can be found in Section \ref{sec:proofs}.

\subsection{Constraint satisfaction}\label{subsec:constraint-satisfaction}
First, we focus on the constraint satisfaction guarantees. We refer to constraint satisfaction as \emph{safety}, which is common in the CBF literature.

\begin{proposition}[Safety]\label{prop:safety}
Let Assumption \ref{assum:differentiability, regularity and relative degree} hold.
If the QP in \eqref{eq:sgf_QP} is feasible everywhere in $\bigcap\limits_{i=0}^{r}\Xc_{i}$ and~\eqref{eq:interconnected-system} has a unique solution
for any initial condition in $\bigcap\limits_{i=0}^{r}\Xc_{i}$, then $\bigcap\limits_{i=0}^{r}\Xc_{i}$ is forward invariant.
\end{proposition} 

In the next section, we provide conditions on feasibility of the QP in \eqref{eq:sgf_QP} and existence and uniqueness of solutions of \eqref{eq:interconnected-system}. As it becomes clear from the proof of Prop. \ref{prop:safety}, all sets $\bigcap_{i=k}^r \{(x,u):\text{ }b(u)\geq 0,h_i(x,u)\geq 0\}$ for all $k=r,r-1,\dots,0$ are forward invariant. As a consequence of Prop. \ref{prop:safety}, if the initial condition lies in $\bigcap\limits_{i=0}^{r}\Xc_{i} \subseteq \mathcal{X}$, the proposed controller enforces both input and state constraints at all times, i.e., $b(\upsilon(t))\geq0$ and $h(\xi(t))\geq 0$, for all $t$ (recall from Remark~\ref{rem:tuning beta} that, under compactness assumptions, sufficiently large values of $\beta$ make  $\bigcap\limits_{i=0}^{r}\Xc_{i} \approx \mathcal{X}$).
Importantly, all feasible points of the optimization problem~\eqref{eq:steady-state_opti_problem} lie in $\bigcap\limits_{i=0}^r\mathcal{X}_i$, as the following result shows.

\begin{proposition}\longthmtitle{Feasible points of~\eqref{eq:steady-state_opti_problem} lie in $\bigcap\limits_{i=0}^r\mathcal{X}_i$}\label{prop:feasible-points-lie-in-bigcap}
    Let $(x_*,u_*)$ be a feasible point of~\eqref{eq:steady-state_opti_problem}. Then, $(x_*,u_*)\in\bigcap\limits_{i=0}^r \mathcal{X}_i$.
\end{proposition}

\subsection{Well-posedness: feasibility and Lipschitz continuity}\label{subsec:well-posedness}
We now provide conditions that guarantee feasibility of the QP \eqref{eq:sgf_QP} and existence and uniqueness
of solutions of \eqref{eq:interconnected-system} over $\Xc_r$. Notice that this suffices to satisfy the assumptions of Prop. \ref{prop:safety}, as $\bigcap\limits_{i=0}^r \mathcal{X}_i\subseteq\Xc_r$.

\begin{proposition}[Feasibility]\label{prop:feasibility}
    Let Assumption \ref{assum:differentiability, regularity and relative degree} hold. Assume that $\mathcal{X}_r$ is compact and that, for any $(x,u) \in \mathcal{X}_r$, there exists $q_{(x,u)}\in\real^m$ such that:
    \begin{itemize}
        \item if $b(u)=0$, then $\nabla^\top b(u) q_{(x,u)} > 0$;
        %
        %
        \item if $h_r(x,u)=0$, then $\frac{\partial h_r}{\partial u}(x,u)q_{(x,u)} + \frac{\partial h_r}{\partial x}(x,u)f(x,u)>0$.
    \end{itemize}
    Then, there exist $\alpha_f, \gamma_f > 0$
    such that, for all $\alpha > \alpha_f$, $\gamma > \gamma_f$, the QP in \eqref{eq:sgf_QP} is feasible for any $(x,u)\in\mathcal{X}_r$.
\end{proposition}

Regarding the hypotheses of Prop. \ref{prop:feasibility}, note first that $\Xc_r$ is compact if  the constraint set $\Xc$ is compact. The result says that sufficient regularity of the set $\{(x,u):b(u)\geq 0, h_r(x,u)\geq 0\}$ guarantees feasibility of QP \eqref{eq:sgf_QP}. To illustrate this, we compare the conditions in Prop. \ref{prop:feasibility} with the standard MFCQ conditions, on the set $\{(x,u): b(u)\geq 0, h_r(x,u)\geq 0\}$. We observe one difference, namely the presence of the term $\frac{\partial h_r}{\partial x}(x,u)f(x,u)$ in the second condition: on the boundary $h_r(x,u)=0$, whenever the plant dynamics $f$ are pushing $h_r(x,u)$ to decrease below 0, the controller's dynamics should be able to counteract that effect. This naturally does not exist in standard MFCQ, as in the case of solving optimization problems, there is no external dynamics (the plant in our case) driving the optimization algorithm outside of the constraint set. Thus, the conditions in Prop. \ref{prop:feasibility} can be seen as MFCQ with the additional assumption that the plant dynamics can be counteracted by the controller in keeping $h_r$ positive.

\begin{proposition}[Lipschitz continuity]\label{prop:Lipschitzness}
Under the assumptions of Prop.~\ref{prop:feasibility}, further assume that,
    for every $(x,u)\in\mathcal{X}_r$, CRCQ holds for~\eqref{eq:sgf_QP}. Let $\alpha_f, \gamma_f$ as in Proposition~\ref{prop:feasibility}.
    Then, by taking $\alpha > \alpha_f$, $\gamma > \gamma_f$, 
    $g_{\epsilon,\alpha,\gamma}$ is locally Lipschitz for all $(x,u)\in\mathcal{X}_r$.
\end{proposition}

The conditions in Prop.~\ref{prop:Lipschitzness} guarantee that~\eqref{eq:interconnected-system} has a unique solution for every initial condition in~$\Xc_r$. For alternative conditions on the optimization problem that can be invoked to guarantee local Lipschitzness,  we refer to~\cite{PM-AA-JC:24-ejc}.

\begin{remark}\longthmtitle{Ensuring CRCQ}
    It can be shown that CRCQ holds for~\eqref{eq:sgf_QP} if Assumption \ref{assum:differentiability, regularity and relative degree} holds, $\mathcal{X}_r$ is compact and there exists $c<1$ such that,  if $b(u) = 0$ and $h_r(x,u) = 0$, then $|\nabla^\top b(u)\frac{\partial h_r}{\partial u}(x,u)|\leq c\|\nabla b(u)\|\cdot\|\frac{\partial h_r}{\partial u}(x,u)\|$, i.e., $\nabla b(u)$ and $\frac{\partial h_r}{\partial u}(x,u)$ are linearly independent in a non-vanishing fashion.
\end{remark}

\subsection{Critical points, equilibria and convergence to optimizers}\label{subsec:equilibria-convergence}
Next, we unveil the relationship between the equilibria of the closed loop \eqref{eq:interconnected-system} and critical points of \eqref{eq:steady-state_opti_problem}.

\begin{proposition}[Critical points of \eqref{eq:steady-state_opti_problem} - equilibria of \eqref{eq:interconnected-system}]\label{prop:critical_points_equilibria}
    Assume that either MFCQ or CRCQ holds at $(x_*,u_*)\in\Xc$ for optimization problem \eqref{eq:steady-state_opti_problem} and QP \eqref{eq:sgf_QP}\footnote{As QP \eqref{eq:sgf_QP} has linear constraints, satisfaction of MFCQ or CRCQ is independent of the value of the decision variable $q$.}, and that at least one of the following holds:
    \begin{enumerate}
        \item $h(x_*) >0$;
        \item $\nabla h(x_*)$ is a right eigenvector of $\frac{\partial f}{\partial x}(x_*,u_*)$ with eigenvalue $e$ and $e^r<0$;
        \item $(\frac{\partial w}{\partial u}(u_*))^\top$ is a left eigenvector of $\frac{\partial f}{\partial x}(x_*,u_*)$ with eigenvalue $e$ and $e^r<0$.
    \end{enumerate}
    Then $(x_*,u_*)$ is a critical point of \eqref{eq:steady-state_opti_problem} if and only if it is an equilibrium of \eqref{eq:interconnected-system}.
\end{proposition}

According to condition 1 in Prop. \ref{prop:critical_points_equilibria}, all critical points in the interior of the state constraint set $\{x:\text{ }h(x)\geq0\}$ are equilibria of \eqref{eq:interconnected-system}, and vice-versa. In practice, this is not too restrictive, as often it is undesirable to operate exactly on the safety boundary $h(x)=0$. Nonetheless, even if this is not the case, we explain in Remark~\ref{rem:perturbing-optimization} how one can modify the controller dynamics \eqref{eq:sgf_QP} so that, given a global optimum of \eqref{eq:steady-state_opti_problem} (even on the boundary), the closed-loop admits a corresponding equilibrium in the interior that is of arbitrarily small suboptimality. 

Equivalence between critical points and equilibria on the boundary is guaranteed if condition 2 or 3 in Prop. \ref{prop:critical_points_equilibria} holds. These conditions can be interpreted as follows. In the steady state, controller \eqref{eq:sgf_QP} can be seen as the SGF solving the optimization problem $\min_{x,u}\Phi(x)$ s.t. $x=w(u), \ b(u)\geq0, \ h_r(x,u)\geq 0$. Either of these conditions enforces that controller \eqref{eq:sgf_QP} is also the SGF of $\min_{u}\Phi(w(u))$ s.t. $b(u)\geq0, \ h(w(u))\geq 0$, which is the original problem \eqref{eq:steady-state_opti_problem}. If these conditions are not satisfied, however, the closed-loop dynamics~\eqref{eq:interconnected-system} might have equilibria on the boundary which are not critical points. Addressing this is left for future work.

Finally, we establish local exponential stability of local minimizers of \eqref{eq:steady-state_opti_problem} under the closed-loop dynamics \eqref{eq:interconnected-system} if condition 1 holds.

\begin{theorem}[Convergence to optima]\label{thm:convergence}
    Let Assumptions \ref{assum:stability} and \ref{assum:differentiability, regularity and relative degree} hold. Let $(x_*,u_*)$ be a local minimizer of \eqref{eq:steady-state_opti_problem} with $h(x_*)>0$ and $\Nc$ a neighborhood where $(x_*,u_*)$ is the only critical point. Assume that $u\mapsto \Phi(w(u))$ has compact level sets. Assume the QP \eqref{eq:sgf_QP} is feasible everywhere in $\Nc$, that MFCQ or CRCQ holds in $\Nc$ for \eqref{eq:steady-state_opti_problem} and~\eqref{eq:sgf_QP}, and that unique solutions of \eqref{eq:interconnected-system} exist for any initial condition in $\Nc$. Then, there exist $\epsilon_*,\gamma_* >0$, such that, for any $\epsilon \in (0,\epsilon_*)$ and $\gamma \geq \gamma_*$, the point $(x_*,u_*)$ is locally asymptotically stable relative\footnote{Asymptotic stability \emph{relative to a set $S$} means that the region of attraction is contained in $S$.} to~$\Xc_r$ for the closed-loop dynamics~\eqref{eq:interconnected-system}.
\end{theorem}

\begin{remark}\longthmtitle{Regularization to obtain interior equilibria}\label{rem:perturbing-optimization}
If condition 1 in Prop. \ref{prop:critical_points_equilibria} is not satisfied, i.e., $h(x_*)=0$, one can do the following. 
    Given $p,\varepsilon >0$, consider the following variation of~\eqref{eq:steady-state_opti_problem}:
    \begin{align}\label{eq:perturbed-optimization-problem}
        \notag
        \min\limits_{x,u}& \ \Phi(x) + p (\varepsilon-h(x)^2) \\
        \mathrm{s.t.}& \ x=w(u), \ b(u) \geq 0, \ h(x) \geq 0.
    \end{align}
    It can be shown that, for any threshold $\delta>0$, there exist $p,\varepsilon>0$ such that the optimizer $(x',u')$ of~\eqref{eq:perturbed-optimization-problem} satisfies $h(x') > 0$ and $|\Phi( x')-\Phi(x_*)| < \delta$. Thus,  using the objective function of \eqref{eq:perturbed-optimization-problem} instead of $\Phi$ in the controller \eqref{eq:sgf_QP}, we can invoke Prop. \ref{prop:critical_points_equilibria} to guarantee that the corresponding closed-loop dynamics admit an equilibrium $(x',u')$ of arbitrarily small suboptimality wrt $(x_*,u_*)$. Further, under the assumptions of Thm. \ref{thm:convergence}, $(x',u')$ can be made locally asymptotically stable relative to $\Xc_r$.
\end{remark}

\section{Numerical Example}\label{sec:example}
Consider the control system
\begin{equation*}
    \begin{aligned}
        &\dot{\xi} = A\xi + B\upsilon,
        \text{ } A = \begin{pmatrix}
            -1.6 &-0.1\\ -1.0& -0.8
        \end{pmatrix}, \text{ } B= \begin{pmatrix}
            1 &0\\ 0 &1
        \end{pmatrix}
    \end{aligned}
\end{equation*}
and observe that it satisfies Assumption \ref{assum:stability}. The steady-state map is $w(u) = -A^{-1}Bu$. We consider a feedback-optimization problem \eqref{eq:steady-state_opti_problem} with the following data:
\begin{equation}\label{eq:example_data}
\begin{aligned}
    &\Phi(x) := \|x-(1.775,0.9)\|^2,\ b(u) := 16-\|u\|^2,\ h(x) := 1 - \frac{(x_1-0.2)^2}{4} - (x_2-0.3)^2, 
\end{aligned}
\end{equation}
which amounts to a strictly convex problem with ellipsoidal state constraints. Assumption \ref{assum:differentiability, regularity and relative degree} is satisfied with $r=1$. The unique global optimizer $(x_*,u_*) = (1.775,0.9)$ lies in the interior of $\Xc$, thus condition 1 in Prop. \ref{prop:critical_points_equilibria} holds, and the conclusions of Thm. \ref{thm:convergence} are valid. 

\begin{figure*}[h!]
\centering
\begin{minipage}[b]{.47\textwidth}
\includegraphics[width=0.99\linewidth]{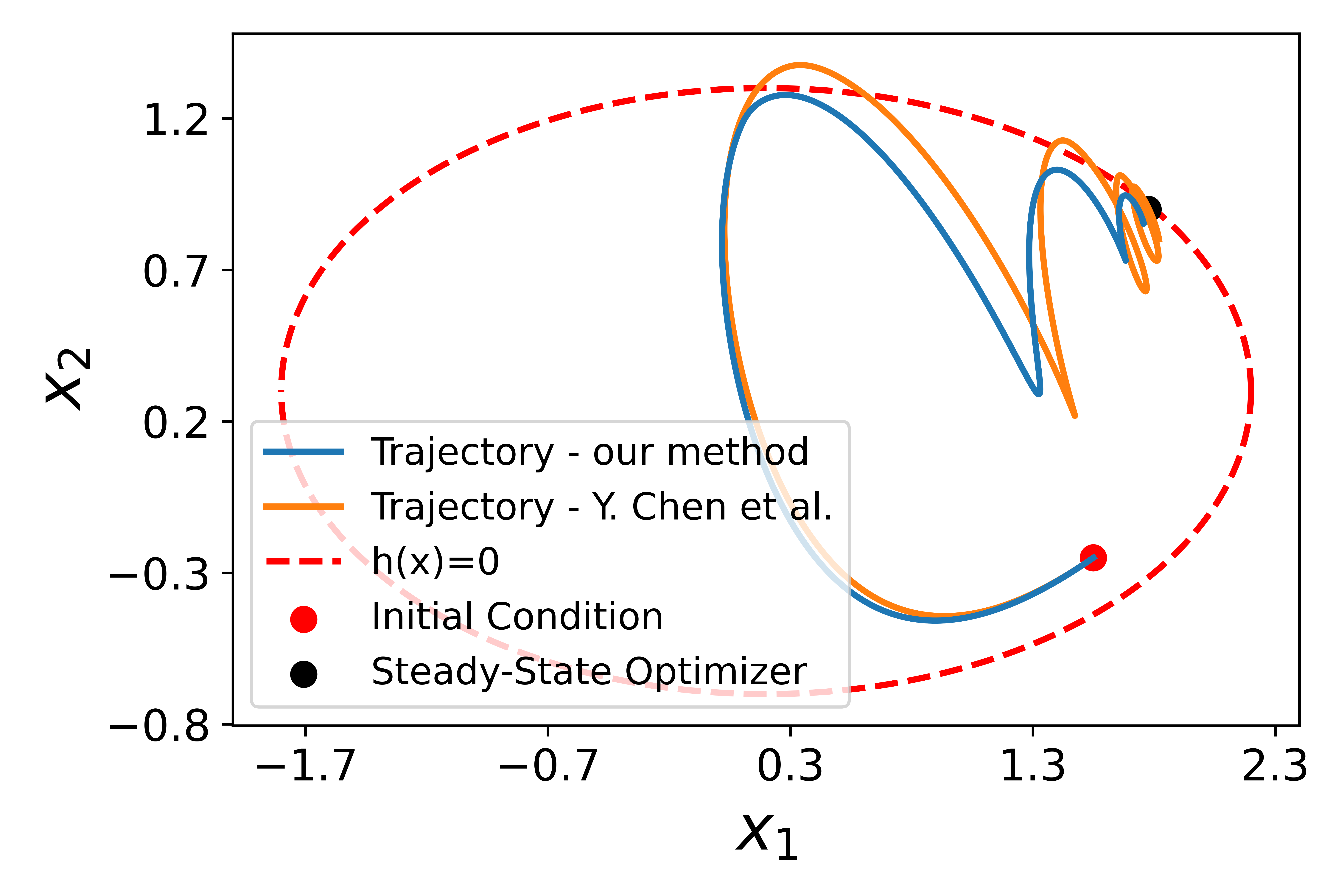}
\caption{Comparison of trajectories of~\eqref{eq:interconnected-system} and the method in~\cite{YC-LC-JC-EDA:23-csl}.}
    \label{fig:1traj}
\end{minipage}\qquad
\begin{minipage}[b]{.47\textwidth}
\includegraphics[width=0.99\linewidth]{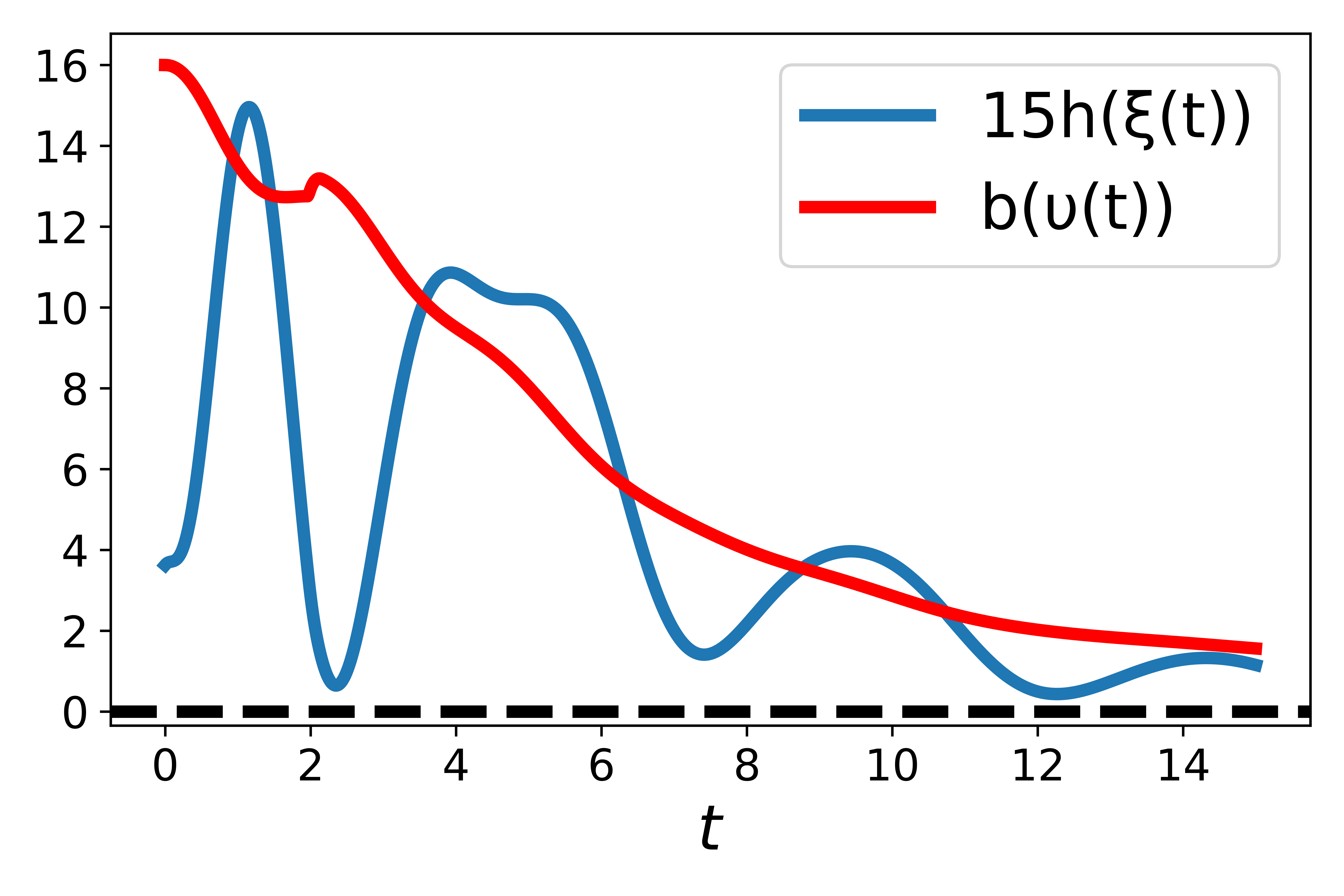}
\caption{Values of $h$ (scaled) and $b$ along a trajectory of~\eqref{eq:interconnected-system}.}
    \label{fig:h_and_b}
\end{minipage}
\end{figure*}

\begin{figure}[h!]
    \centering
    \includegraphics[width=0.47\linewidth]{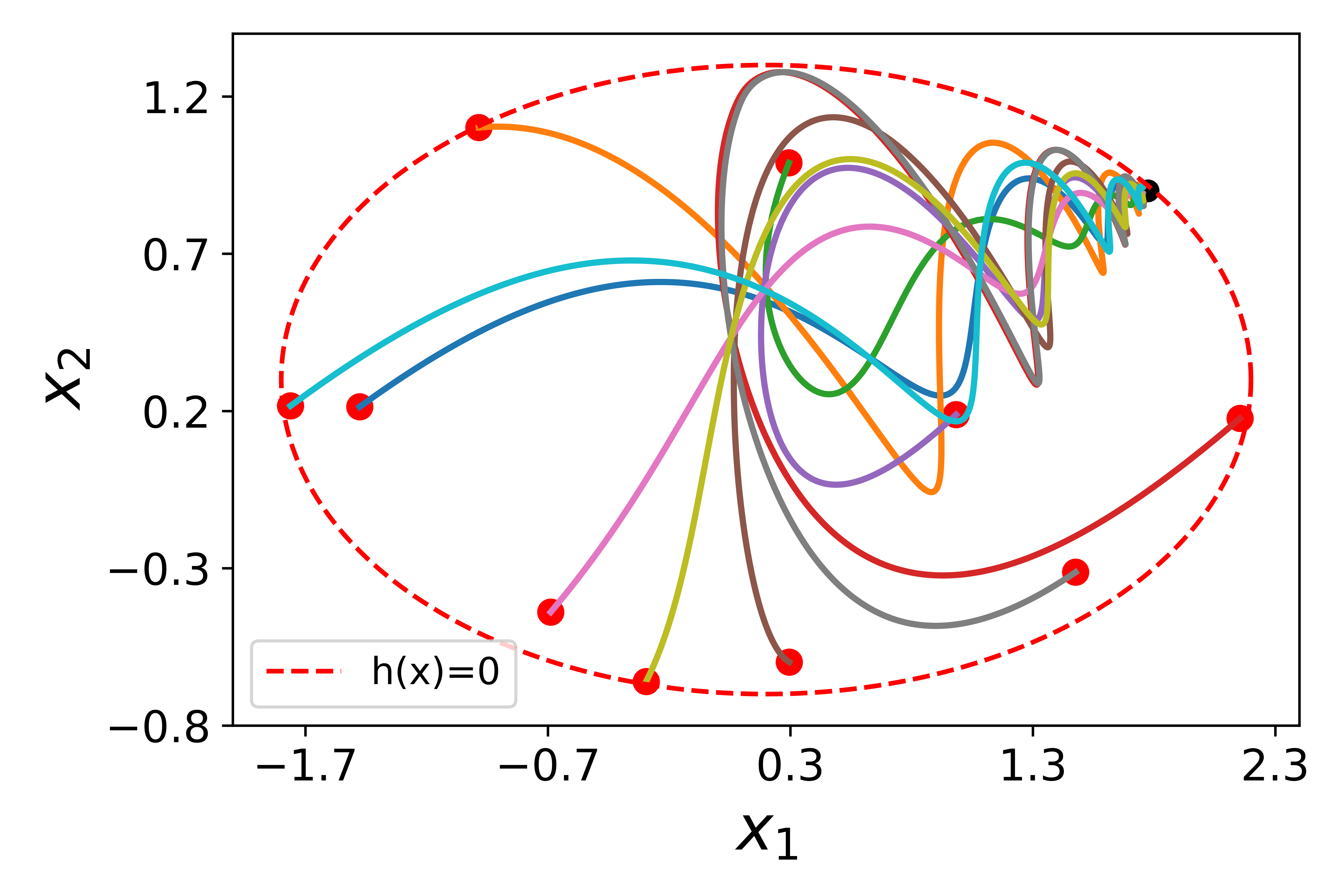}
    \caption{Trajectories of~\eqref{eq:interconnected-system} for different initial conditions.}
    \label{fig:multiple_traj}
\end{figure}

Figure \ref{fig:1traj} compares our approach with the one proposed in~\cite{YC-LC-JC-EDA:23-csl}, that handles only input constraints. Both approaches converge to the steady-state optimizer. However, as expected, the latter leads to transient state-constraint violations, while the design proposed here enforces them at all times; as Figure \ref{fig:h_and_b} shows, our method, accounting for dynamic plants and employing high-order CBFs, enforces $h(\xi(t))\geq 0$ and $b(\upsilon(t))\geq 0$ for all $t$. Finally, in Figure \ref{fig:multiple_traj}, we simulate our approach for several initial conditions. Even though Thm. \ref{thm:convergence} guarantees only local stability, it seems that, in this example, stability is global wrt the constraint set. This may be related to the strong convexity of feedback optimization problem \eqref{eq:steady-state_opti_problem}-\eqref{eq:example_data} and will be subject of future research.

\section{Conclusions}
We have presented a feedback-optimization method that handles both input and state constraints at all times. We have provided results on the proposed dynamics' well-posedness, constraint satisfaction, equivalence between equilibria and critical points, and local stability of local optima. Aspects in need of further research include: a) there might exist spurious equilibria on the state-constraint set's boundary that are not critical points, b) the stability guarantee of optima is only local, c) the system model is used to enforce state constraints. Future research will focus on (partially) model-free cases and on guaranteeing global convergence to optima or critical points, thereby excluding limit cycles and convergence to undesired equilibria.

\section{Technical proofs}\label{sec:proofs}
\begin{proof}[\textbf{Proof of Prop. \ref{prop:safety}}]
    If, for some $t'>0$, we have $b(\upsilon(t'))\geq0$, then the first constraint in \eqref{eq:sgf_QP} enforces that $b(\upsilon(t))\geq 0$ for all $t\geq t'$ (through standard CBF arguments, see e.g.~\cite{xiao2021high}). Similarly, if for some $t'>0$, we have $h_r(\xi(t'),\upsilon(t'))\geq0$, the second constraint enforces that $h_r(\xi(t),\upsilon(t))\geq 0$ for all $t\geq t'$, i.e., the set $\{(x,u):h_r(x,u)\geq 0\}$ is forward invariant. Now, if for $t'>0$, we have $h_{r-1}(\xi(t'),\upsilon(t'))\geq 0$, again through standard CBF arguments, we have that $h_{r-1}(\xi(t),\upsilon(t))\geq 0$ for all $t\geq t'$ if $h_r(\xi(t),\upsilon(t))\geq 0$ for all $t\geq t'$, due to the definition of $h_r$. Thus, the set $\{(x,u):h_r(x,u)\geq 0, \ h_{r-1}(x,u)\geq 0\}$ is forward invariant. Using this argument recursively, we conclude that $\bigcap\limits_{i=0}^r\mathcal{X}_i$ is forward invariant.
\end{proof}

\begin{proof}[\textbf{Proof of Prop. \ref{prop:feasible-points-lie-in-bigcap}}]
    Since $(x_*,u_*)$ is feasible for~\eqref{eq:steady-state_opti_problem}, $h(x_*) \geq 0$ and $b(u_*) \ge 0$. Therefore, $(x_*,u_*)\in\mathcal{X}_0$.
    Now, note that $h_1(x_*,u_*) = \nabla h(x_*,u_*)^T f(x_*,u_*) + \beta h(x_*)$, and since $f(x_*,u_*)=f(w(u_*),u_*)=0$ and $\beta > 0$, we obtain $h_1(x_*,u_*) \geq 0$.
    In fact, for all $i\in[r]$, $h_i(x_*,u_*) = \beta h_{i-1}(x_*,u_*)$, implying that $h_i(x_*,u_*) \geq 0$ for all $i\in[r]$ and hence $(x_*,u_*)\in\bigcap\limits_{i=0}^r \mathcal{X}_i$.
\end{proof}

\begin{proof}[\textbf{Proof of Prop. \ref{prop:feasibility}}]
    For each $(x,u)\in\partial\Xc_r$, let $\Nc_{(x,u)}$ be an open neighborhood of $(x,u)$ such that for every $(\bar{x},\bar{u})\in\Nc_{(x,u)}$
    \begin{itemize}
        \item if $b(u)=0$, then $\nabla^\top b(\bar{u}) q_{(x,u)} > 0$;
        \item if $h_r(x,u)=0$, then $\frac{\partial h_r}{\partial u}(\bar{x},\bar{u})q_{(x,u)} + \frac{\partial h_r}{\partial x}(\bar{x},\bar{u})f(\bar{x},\bar{u})>0$.
    \end{itemize}
    Note that such a neighborhood $\Nc_{(x,u)}$ exists by continuity of $\nabla b$, $\frac{\partial h_r}{\partial u}$ and $(x,u)\mapsto \frac{\partial h_r}{\partial x}(x,u)f(x,u)$.
     
    Now, let $S_1 := \setdef{(x,u)\in\Xc_r}{b(u) = 0}$, $S_2 := \setdef{(x,u)\in\Xc_r}{h_r(x,u) = 0}$, $S_3 :=S_1\cap S_2$, and
    consider the following three sets:
    \begin{align*}
        &\Nc_1 = \hspace{-2mm}\bigcup_{ (x,u)\in S_1 } \hspace{-2mm}\Nc_{(x,u)}, \ \Nc_2 = \hspace{-2mm}\bigcup_{ (x,u)\in S_2 }\hspace{-2mm} \Nc_{(x,u)}, \ \Nc_3 = \hspace{-2mm}\bigcup_{ (x,u)\in S_3 } \hspace{-2mm}\Nc_{(x,u)}.
    \end{align*}
    %
    %
    Since $\Nc_1$ is a cover of $S_1$, we can extract a finite sub-cover $\Nc_{f,1} = \bigcup_{i=1}^{N_1} \Nc_{(x_i,u_i)}$ of $S_1$ with $b(u_i) = 0$ for all $i\in [N_1]$ and $N_1 \in \mathbb{Z}_{>0}$.
    Similarly, we can extract finite sub-covers 
    $\Nc_{f,2} = \bigcup_{i=1}^{N_2} \Nc_{(x_i,u_i)}$ of $S_2$, with $h_r(x_i,u_i) = 0$ for all $i\in[N_2]$ and $N_2 \in \mathbb{Z}_{>0}$, and 
    $\Nc_{f,3} = \bigcup_{i=1}^{N_3} \Nc_{(x_i,u_i)}$ of $S_3$, with $b(u_i) = h_r(x_i,u_i) = 0$ for all $i\in[N_3]$ and $N_3 \in \mathbb{Z}_{>0}$.

    Since $b$ is continuous and $S_1$ is compact, there exists $c_1>0$ such that, if $b(u) \in [0, c_1]$, then $(x,u)\in\Nc_{f,1}$. Indeed, if this was not the case, there would exist a sequence of positive numbers $\{ \bar{c}_i \}_{i\in\mathbb{Z}_{>0}}$ with $\lim\limits_{i\to\infty} \bar{c}_i = 0$ and points $(\bar{x}_i,\bar{u}_i)\in \Xc_r$ such that $b(\bar{u}_i) \in [0, \bar{c}_i]$ and $(\bar{x}_i,\bar{u}_i)\notin\Nc_{f,1}$. As $\Xc_r\setminus\Nc_{f,1}$ is compact, a subsequence of these points converges to a point $(x_*,u_*)\in \Xc_r\setminus\Nc_{f,1}$ with $b(u_*)=0$; this is a contradiction, as $b(u_*)=0$ implies $(x_*,u_*)\in S_1\subseteq \Nc_{f,1}$. 
    Similarly, there exists $c_2>0$ such that, if $h_r(x,u)\in[0,c_2]$, then $(x,u)\in\Nc_{f,2}$; and $c_3>0$ such that, if $b(u) \in [0,c_3]$ and $h_r(x,u)\in[0,c_3]$, then $(x,u)\in\Nc_{f,3}$.

    Now we consider four cases:
    \begin{enumerate}
        \item Let $Q_3 = \setdef{(x,u)\in\Xc_r}{h_r(x,u)\in[0,c_3], b(u)\in[0,c_3]}$. Since $Q_3 \subset \Nc_{f,3}$, 
        then, for each $(x,u)\in Q_3$,
        there exists $i\in[N_3]$ such that 
        \begin{align*}
            &\nabla b(u)^\top q_{(x_i,u_i)} + \alpha b(u) > 0, \\
            &\frac{\partial h_r}{\partial u}(x,u) q_{(x_i,u_i)} + \frac{\partial h_r}{\partial x}(x,u)f(x,u) + \gamma h_r(x,u) > 0,
        \end{align*}
        for all $\alpha\geq 0$, $\gamma\geq0$. Thus, QP \eqref{eq:sgf_QP} is feasible in $Q_3$, for all $\alpha\geq 0$, $\gamma\geq0$.
        
        \item Let $Q_2 = \setdef{(x,u)\in\Xc_r}{ h_r(x,u)\in[0,c_2], b(u)>c_3 }$. Since $Q_2 \subset \Nc_{f,2}$, for each $(x,u)\in Q_2$, we know that there exists $i\in[N_2]$ such that
        \begin{equation*}
            \frac{\partial h_r}{\partial u}(x,u) q_{(x_i,u_i)} + \frac{\partial h_r}{\partial x}(x,u)f(x,u) + \gamma h_r(x,u) > 0,
        \end{equation*}
        for any $\gamma\geq0$. Define 
        \begin{align*}
            \alpha_{f,1} = \frac{ \max\limits_{ (x,u)\in \mathcal{X}_r, i\in\{1,\hdots,N_2 \} } |\nabla b(u)^\top q_{(x_i,u_i)}| }{c_3}.
        \end{align*}
        Since $\Xc_r$ is compact, $\nabla b$ is continuous, and $N_2$ is finite, $\alpha_{f,1}$ is finite. Therefore $\nabla b(u)^\top q_{(x_i,u_i)} + \alpha b(u) \geq 0$ for any $\alpha > \alpha_{f,1}$. Thus, QP \eqref{eq:sgf_QP} is feasible in $Q_2$, for all $\alpha > \alpha_{f,1}$, $\gamma\geq0$.

        \item Let $Q_1 = \setdef{(x,u)\in \Xc_r}{b(u)\in[0,c_1], h_r(x,u)>c_3 }$. Since $Q_1 \subset \Nc_{f,1}$, for each $(x,u)\in Q_1$, we know that there exists $i\in[N_1]$ such that $\nabla b(u)^\top q_{(x_i,u_i)} + \alpha b(u) > 0$ for all $\alpha\geq 0$. Define
        \begin{align*}
            \gamma_{f,1} = \frac{1}{c_3} \max\limits_{ (x,u)\in \mathcal{X}_r, i\in\{1,\hdots,N_1 \} } \Big| \frac{\partial h_r}{\partial u}(x,u) q_{(x_i,u_i)} + \frac{\partial h_r}{\partial x}(x,u)f(x,u) \Big|.
        \end{align*}
        Since $\Xc_r$ is compact, $\frac{\partial h_r}{\partial u}$ and $(x,u)\mapsto\frac{\partial h_r}{\partial x}(x,u)f(x,u)$ are continuous, and $N_1$ is finite, $\gamma_{f,1}$ is finite. Therefore $\frac{\partial h_r}{\partial u}(x,u) q_{(x_i,u_i)} + \frac{\partial h_r}{\partial x}(x,u)f(x,u) + \gamma h_r(x,u) > 0$, for all $\gamma > \gamma_{f,1}$. Thus, QP \eqref{eq:sgf_QP} is feasible in $Q_1$, for all $\alpha\geq0$, $\gamma > \gamma_{f,1}$.
        
        \item Let $Q_4 = \mathcal{X}_r \setminus ( Q_1 \cup Q_2 \cup Q_3 )$.
        Note that for all $(x,u) \in Q_4$, we have 
        $b(u) > \min\{ c_1, c_3 \} > 0$ and 
        $h_r(x,u) > \min\{ c_2, c_3 \} > 0$.
        %
        %
        Define 
        \begin{align*}
            \gamma_{f,2} = \frac{ \max\limits_{ (x,u)\in (\mathcal{X}\times\mathcal{U})\cap\mathcal{X}_r } | \frac{\partial h_r}{\partial x}(x,u)f(x,u) | }{ \min\{ c_2, c_3 \} }.
        \end{align*}
        Since $\mathcal{X}_r$ is compact, and $\frac{\partial h_r}{\partial x}$ and $f$ are continuous, $\gamma_{f,2}$ is finite. Now, for all $(x,u)\in Q_4$ we have 
        \begin{align*}
            &\nabla b(u)^\top \mathbf{0}_m + \alpha b(u) > 0, \\
            &\frac{\partial h_r}{\partial u}(x,u)^\top \mathbf{0}_m + \frac{\partial h_r}{\partial x}(x,u)f(x,u) + \gamma h_r(x,u) > 0.
        \end{align*}
        for all $\alpha \geq 0$ and $\gamma \geq \gamma_{f,2}$. Thus, QP \eqref{eq:sgf_QP} is feasible in $Q_4$, for all $\alpha > 0$, $\gamma > \gamma_{f,2}$.
    \end{enumerate}
    The result follows by observing that $\Xc_r = Q_1\cup Q_2 \cup Q_3 \cup Q_4$, and taking $\alpha_f = \alpha_{f,1}$, $\gamma_f = \max\{ \gamma_{f,1}, \gamma_{f,2} \}$.
\end{proof}

\begin{proof}[\textbf{Proof of Prop. \ref{prop:Lipschitzness}}]
    By following the argument used in Prop.~\ref{prop:feasibility}, 
    by taking $\alpha > \alpha_f$, $\gamma > \gamma_f$,~\eqref{eq:sgf_QP} is strictly feasible (i.e., for any $(x,u)\in\mathcal{X}_r$ there exists $q\in\real^m$ satisfying each of the constraints strictly). Thus, Slater's condition holds for~\eqref{eq:sgf_QP}. Since~\eqref{eq:sgf_QP} is a strongly convex QP, by~\cite[Prop. 5.39]{NA-AE-MP:20}, we have that MFCQ holds for~\eqref{eq:sgf_QP} at the optimizer. The result follows from~\cite[Thm. 3.6]{JL:95}.
\end{proof}

\begin{proof}[\textbf{Proof of Prop. \ref{prop:critical_points_equilibria}}]
        \textbf{Proof of $\implies$.} First, because $x_*=w(u_*)$, we have that $f(x_*,u_*)=f(w(u_*),u_*)=0$, from Assumption \ref{assum:differentiability, regularity and relative degree}. Due to the assumption on MFCQ or CRCQ, there exist $\lambda_h,\lambda_b\geq 0$ and $\mu\in\real^n$ such that the KKT conditions for \eqref{eq:steady-state_opti_problem} are satisfied:
    \begin{align}\label{eq:equiibria_proof_kkt_feedback_opti}
        &\begin{bmatrix}
            \nabla \Phi(x_*) \\ 0 
        \end{bmatrix}  - \lambda_h \begin{bmatrix}
            \nabla h(x_*) \\ 0 
        \end{bmatrix} - \lambda_b \begin{bmatrix}
            0 \\ \nabla b(u_*) 
        \end{bmatrix} + 
        \begin{bmatrix}
            \mu\\ -\dfrac{\partial w}{\partial u}(u_*)^\top\mu
        \end{bmatrix} = 0,\\
        &0\leq \lambda_h \perp h(x_*)\geq0, \quad 0\leq \lambda_b \perp b(u_*)\geq0, \quad x_* = w(u_*) \notag
    \end{align}
    Further, at $(x_*,u_*)$, we have:
    \begin{equation*}
        g_{\epsilon,\alpha,\gamma}(x_*,u_*)=\left\{\begin{aligned}&\argmin_{q}\text{ }\frac{1}{2}\|q +\frac{\partial w}{\partial u}(u_*)^\top \nabla\Phi(x_*)\|^2\\
            &\quad\mathrm{s.t.:} \text{ } \nabla^\top b(u_*)\cdot q +\alpha b(u_*)\geq 0\\
            &\quad\quad\text{ }\begin{aligned} \underbrace{\frac{\partial h_r}{\partial u}(x_*,u_*)}_{\frac{\partial \lie^r_{f}h}{\partial u}(x_*,u_*)}\cdot q+\gamma h_{r}(x_*,u_*)\geq 0&\end{aligned}\end{aligned}\right.
    \end{equation*}
    Now, observe that $ \tilde{f}(u)=f(w(u),u)=0$. Thus $\dfrac{\partial \tilde{f}}{\partial u}(u)=0\implies \dfrac{\partial f}{\partial u}(w(u),u)=-\dfrac{\partial f}{\partial x}(w(u),u)\dfrac{\partial w}{\partial u}(u)$. By employing this relationship, we obtain:
    \begin{equation}\label{eq:using_partial_of_w}
        \Big(\dfrac{\partial}{\partial u}\lie^r_{f}h(x_*,u_*)\Big)^\top = -\dfrac{\partial w}{\partial u}(u_*)^\top\Big(\dfrac{\partial f}{\partial x}(x_*,u_*)\Big)^r\nabla h(x_*)
    \end{equation}
    
    The KKT conditions for the controller's QP are as follows:\begin{equation}\label{eq:equiibria_proof_kkt_controller}
    \begin{aligned}
        &q_* +\Big(\lambda_h\frac{\partial w}{\partial u}(u_*)^\top \nabla h(x_*) +\lambda_b\nabla b(u_*) \Big)- \lambda'_b\nabla b(u_*) +\lambda'_h\dfrac{\partial w}{\partial u}(u_*)^\top\Big(\dfrac{\partial f}{\partial x}(x_*,u_*)\Big)^r\nabla h(x_*)= 0\\
        &0\leq\lambda'_b \perp \nabla^\top b(u_*)\cdot q_* +\alpha b(u_*)\geq 0\\ 
        &0\leq\lambda'_h \perp \frac{\partial \lie^r_{f}h}{\partial u}(x_*,u_*)\cdot q_*+\beta^{r}\gamma h(x_*)\geq0\\
    \end{aligned}
    \end{equation}
    where $\lambda'_h,\lambda'_b$ are Lagrange multipliers; we replaced $\frac{\partial w}{\partial u}(u_*)^\top \nabla\Phi(x_*)$ by  $\lambda_h\frac{\partial w}{\partial u}(u_*)^\top \nabla h(x_*) +\lambda_b\nabla b(u_*)$, employing \eqref{eq:equiibria_proof_kkt_feedback_opti}; we used \eqref{eq:using_partial_of_w}; we used that $h_r(x_*,u_*) = \beta h_{r-1}(x_*,u_*) = \beta^2 h_{r-2}(x_*,u_*)=\dots=\beta^r h(x_*)$, employing \eqref{eq:h_i} and $f(x_*,u_*) = 0$. If we can find $\lambda'_h,\lambda'_b \geq 0$ such that \eqref{eq:equiibria_proof_kkt_controller} holds with $q_* = 0$, the statement is proven. 
    \begin{itemize}
        \item If condition (1) holds, i.e., $h(x_*)>0$, then $\lambda_h = 0$. Picking $\lambda'_h = 0$, $\lambda'_b = \lambda_b$, \eqref{eq:equiibria_proof_kkt_controller} is satisfied with $q_* = 0$.
        \item If any of (2) or (3) holds, then picking $\lambda'_h = -\frac{1}{e^r}\lambda_h$, $\lambda'_b = \lambda_b$, we see that \eqref{eq:equiibria_proof_kkt_controller} is satisfied with $q_* = 0$.
    \end{itemize}

    \textbf{Proof of $\impliedby$.} As $(x_*,u_*)$ is an equilibrium of \eqref{eq:interconnected-system}, we have $x_* = w(u_*)$, $f(x_*, u_*)=0$ and $g_{\epsilon,\alpha,\gamma}(x_*, u_*) = 0$. Through the KKT conditions of the controller's QP, and noting that $q_* = g_{\epsilon,\alpha,\gamma}(x_*, u_*) = 0$, we obtain that there exist $\lambda'_h,\lambda'_b\geq 0$ that satisfy the following conditions:
    \begin{equation}\label{eq:equilibria_proof_controller_kkt2}
    \begin{aligned}
        &\frac{\partial w}{\partial u}(u_*)^\top \nabla\Phi(x_*)- \lambda'_b\nabla b(u_*) +\lambda'_h\dfrac{\partial w}{\partial u}(u_*)^\top\Big(\dfrac{\partial f}{\partial x}(x_*,u_*)\Big)^r\nabla h(x_*)= 0\\
        &0\leq\lambda'_b \perp \alpha b(u_*)\geq 0, \quad 0\leq\lambda'_h \perp \beta^{r}\gamma h(x_*)\geq0\\
    \end{aligned}
    \end{equation}
    where we used \eqref{eq:using_partial_of_w} and that $h_r(x_*,u_*) = \beta h_{r-1}(x_*,u_*) = \beta^2 h_{r-2}(x_*,u_*)=\dots=\beta^r h(x_*)$. To prove the statement, we need to find $\lambda_h,\lambda_b\geq 0$ and $\mu\in\real^n$ such that \eqref{eq:equiibria_proof_kkt_feedback_opti} holds. We pick $\mu = -\nabla\Phi(x_*) + \lambda_h \nabla h(x_*)$, and \eqref{eq:equiibria_proof_kkt_feedback_opti} becomes:
    \begin{equation*}
    \begin{aligned}
        &\begin{bmatrix}
            0 \\ 0 
        \end{bmatrix}  - \lambda_h \begin{bmatrix}
            0 \\ 0 
        \end{bmatrix} - \lambda_b \begin{bmatrix}
            0 \\ \nabla b(u_*) 
        \end{bmatrix} +\begin{bmatrix}
            0\\ -\dfrac{\partial w}{\partial u}(u_*)^\top(-\nabla\Phi(x_*) + \lambda_h \nabla h(x_*))
        \end{bmatrix} = 0\\
        &0\leq \lambda_h \perp h(x_*)\geq0, \quad 0\leq \lambda_b \perp b(u_*)\geq0, \quad x_* = w(u_*)
    \end{aligned}
    \end{equation*}
    Employing the first line of \eqref{eq:equilibria_proof_controller_kkt2}, the above become:
    \begin{equation}\label{eq:equiibria_proof_kkt_feedback_opti2}
    \begin{aligned}
        &-\lambda_b\nabla b(u_*) -\dfrac{\partial w}{\partial u}(u_*)^\top\lambda_h \nabla h(x_*) +\lambda'_b\nabla b(u_*) -\lambda'_h\dfrac{\partial w}{\partial u}(u_*)^\top\Big(\dfrac{\partial f}{\partial x}(x_*,u_*)\Big)^r\nabla h(x_*) = 0\\
        &0\leq \lambda_h \perp h(x_*)\geq0, \quad 0\leq \lambda_b \perp b(u_*)\geq0, \quad x_* = w(u_*)
    \end{aligned}
    \end{equation}
    Again, we have the following cases:
    \begin{itemize}
        \item If condition (1) holds, i.e., $h(x_*)>0$, then $\lambda'_h = 0$. Picking $\lambda_h = 0$, $\lambda'_b = \lambda_b$, we see that \eqref{eq:equiibria_proof_kkt_feedback_opti2} holds.
        \item If any of (2) or (3) holds, then picking $\lambda_h = -e^r\lambda'_h$, $\lambda'_b = \lambda_b$, we see that \eqref{eq:equiibria_proof_kkt_feedback_opti2} holds.
    \end{itemize}
\end{proof}

\begin{proof}[\textbf{Proof of Thm. \ref{thm:convergence}}]
    First, from~\cite[Lemma 2.1]{YC-LC-JC-EDA:23-csl}, under Assumption \ref{assum:stability}, there exists $W:\real^n\times\real^m\to\real$, and constants $d_1,d_2,d_3,d_4>0$ such that for all $(x,u)\in\real^n\times\real^m$:
\begin{align*}
        &d_1 \norm{x-w(u)}^2 \leq W(x,u) \leq d_2 \norm{x-w(u)}^2, \\
        &\frac{\partial W(x,u)}{\partial x}f(x,u) \leq -d_3 \norm{x-w(u)}^2, \\
        &
        \norm{\frac{\partial W(x,u)}{\partial u}} \leq d_4 \norm{x-w(u)}.
    \end{align*}
    
    Consider the function $V:\Nc\to\real$ defined by $V(x,u):= \Phi(w(u))-\Phi(w(u_*))+W(x,u)$. Let us show that $V$ is a Lyapunov function. Note that $V(x_*,u_*) = 0$, as $x_*=w(u_*)$, and 
    since
    $\Phi(w(u))>\Phi(w(u_*))$ for all $(x,u)\in\Nc\setminus\{(x_*,u_*)\}$,
    $V(x,u) > 0$ for all $(x,u)\in \Nc \setminus\{ (x_*,u_*) \}$. 
    
    Now, we have:
    \begin{equation*}
       \frac{d}{dt} \Phi(w(\upsilon(t)))=  \nabla^\top\Phi(w(\upsilon(t))) \frac{\partial w}{\partial u}(\upsilon(t)) g_{\epsilon,\alpha,\gamma}(\xi(t),\upsilon(t)),
    \end{equation*}
    where $(\xi(t),\upsilon(t))$ is a solution to \eqref{eq:interconnected-system}, with initial condition $(\xi(0),\upsilon(0))\in \Xc_r$. Notice that $b(\upsilon(t))\geq 0$ and $h_r(\xi(t),\upsilon(t))\geq 0$ for all time, from forward invariance of $\{(x,u): \ b(u)\geq0, \ h_r(x,u) \geq 0\}$ (see proof of Prop. \ref{prop:safety}). In the following, we drop the dependence of $\xi$ and $\upsilon$ on $t$ from notation, for convenience.
    Let $\Gamma\subseteq \Nc\cap\Xc_r$ be a sublevel set of $V$, such that $(x_*,u_*)\in \Gamma$ and $h_r(x,u) > 0$ for all $(x,u)\in\Gamma$.\footnote{Such a sublevel set exists because: a) $V$ has compact sublevel sets due to $\Phi$ having compact sublevel sets and $W$ being radially unbounded; and b) $(x_*,u_*)$ satisfies $h_r(x_*,u_*)=\beta^r h(x_*)>0$, and due to continuity of $h_r$, $\Gamma$ can be picked small enough so that $h(x)>0$ for all $(x,u)\in\Gamma$.} Since $\nabla \Phi$ and $\frac{\partial w}{\partial u}$ are locally Lipschitz, 
    and $\Gamma$ is compact, there exists $L>0$ such that for all $(x,u),(x^\prime,u)\in\Gamma$, we have
    $\|\nabla^\top \Phi(x)\frac{\partial w}{\partial u}(u)-\nabla^\top \Phi(x')\frac{\partial w}{\partial u}(u)\|\leq L\|x-x'\|$. Then, whenever $(\xi(t),\upsilon(t))\in\Gamma$,
    \begin{equation}\label{eq:giannis1}
        \begin{aligned}
            \frac{d}{dt} \Phi(w(\upsilon))=&\\\big(\nabla^\top\Phi(w(\upsilon))\cdot\frac{\partial w}{\partial u}(\upsilon)-\nabla^\top\Phi(\xi)\cdot\frac{\partial w}{\partial u}(\upsilon)\big)g_{\epsilon,\alpha,\gamma}(\xi,\upsilon)+\nabla^\top\Phi(\xi)\cdot\frac{\partial w}{\partial u}(\upsilon)g_{\epsilon,\alpha,\gamma}(\xi,\upsilon)\leq&\\
        L \|\xi-w(\upsilon)\|\|g_{\epsilon,\alpha,\gamma}(\xi,\upsilon)\|
         +\nabla^\top\Phi(\xi)\cdot\frac{\partial w}{\partial u}(\upsilon)g_{\epsilon,\alpha,\gamma}(\xi,\upsilon)\text{ }&
        \end{aligned}
    \end{equation}
    Now, from KKT conditions of~\eqref{eq:sgf_QP} we have:
   \begin{align}\label{eq:giannis_kkt1}
        g_{\epsilon,\alpha,\gamma}(x,u) +  \epsilon\frac{\partial w}{\partial u}(u)^\top\nabla\Phi(x)&-\lambda'_b(x,u)\nabla b(u) - \lambda'_h(x,u) \frac{\partial h_r}{\partial u}(x,u)^\top = 0 
    \end{align}
    where $\lambda'_b,\lambda'_h$ are Lagrange multipliers satisfying:
\begin{align}\label{eq:giannis_kkt2}
        &0\leq\lambda'_b(x,u)\perp \nabla^\top b(u)g_{\epsilon,\alpha,\gamma}(x,u) + \alpha b(u)\geq 0, 
        \\
        &0\leq\lambda'_h(x,u)\perp \frac{\partial h_r}{\partial u}(x,u)g_{\epsilon,\alpha,\gamma}(x,u) + \frac{\partial h_r}{\partial x}(x,u)f(x,u) +\gamma h_r(x,u)\geq 0 \notag
    \end{align}
    For the second term of \eqref{eq:giannis1}, we obtain:
    \begin{equation}\label{eq:giannis2}
    \begin{aligned}
    \nabla^\top\Phi(\xi)\cdot\frac{\partial w}{\partial u}(\upsilon)g_{\epsilon,\alpha,\gamma}(\xi,\upsilon)=&\\ \frac{1}{\epsilon} \Big(-g_{\epsilon,\alpha,\gamma}(\xi,\upsilon) + \lambda'_b(\xi,\upsilon)\nabla b(\upsilon)+  \lambda'_h(\xi,\upsilon)\frac{\partial h_r}{\partial u}(\xi,\upsilon(t))\Big)^\top g_{\epsilon,\alpha,\gamma}(\xi,\upsilon)
    =&\\ - \frac{1}{\epsilon} \|g_{\epsilon,\alpha,\gamma}(\xi,\upsilon)\|^2 - \frac{1}{\epsilon}\lambda'_b(\xi,\upsilon)\alpha b(\upsilon) - 
    \frac{1}{\epsilon}\lambda'_h(\xi,\upsilon)\gamma h_r(\xi,\upsilon)-\frac{1}{\epsilon}\lambda'_h(\xi,\upsilon)\frac{\partial h_r}{\partial x}(\xi,\upsilon)f(\xi,\upsilon), \ &
    \end{aligned} 
    \end{equation}
    whenever $(\xi(t),\upsilon(t))\in\Gamma$, where in the first step we used \eqref{eq:giannis_kkt1} and in the second step we used \eqref{eq:giannis_kkt2}.
    Now, since $h_r(x,u) > 0$ for all $(x,u)\in\Gamma$, and $\Gamma$ is compact,
    there exists $c_0>0$ such that $h_r(x,u)\geq c_0$
    for all $(x,u)\in\Gamma$. Moreover, since $\Gamma$ is compact and $\frac{\partial h_r}{\partial x}$ and $f$ are continuous, there exists $c_1>0$ such that 
    $\frac{\partial h_r}{\partial x}(x,u)f(x,u) \geq -c_1$ for all $(x,u)\in\Gamma$. Hence, by taking $\gamma \geq \gamma_* := \frac{c_1}{c_0}$, 
    we have $\gamma h_r(x,u) + \frac{\partial h_r}{\partial x}(x,u)f(x,u) \geq 0$ for all $(x,u)\in\Gamma$, which from~\eqref{eq:giannis2} means that 
    \begin{align*}
        \nabla^\top\Phi(\xi)\cdot\frac{\partial w}{\partial u}(\upsilon)g_{\epsilon,\alpha,\gamma}(\xi,\upsilon) \leq -\frac{1}{\epsilon}\|g_{\epsilon,\alpha,\gamma}(\xi,\upsilon)\|^2,
    \end{align*}
    whenever $(\xi,\upsilon)\in\Gamma$,
    where we have also used the fact that $b(\upsilon(t)) \geq 0$ for all $t\geq0$.
    Note also that
    \begin{equation*}
        \begin{aligned}
            \frac{d}{dt}W(\xi,\upsilon) \! \leq \! & - \! d_3 \! \norm{\xi-w(\upsilon)}^2 \! + \! d_4 \! \norm{ g_{\epsilon,\alpha,\gamma}(\xi,\upsilon)} \norm{\xi-w(\upsilon)}.
        \end{aligned}
    \end{equation*}
    
    Hence, we have $\frac{d}{dt}V(\xi(t),\upsilon(t)) \leq \zeta(t)^T P_{\epsilon} \zeta(t)$
    whenever $(\xi(t),\upsilon(t))\in\Gamma$, where
    \begin{align*}\zeta(t) = \begin{pmatrix} \norm{g_{\epsilon,\alpha,\gamma}(\xi(t),\upsilon(t))} &\norm{\xi(t)-w(\upsilon(t))} \end{pmatrix}, \ 
        P_{\epsilon} = \begin{pmatrix}
            -\frac{1}{\epsilon} & \frac{ (L+d_4) }{2} \\
            \frac{(L+d_4) }{2} & -d_3
        \end{pmatrix}.
    \end{align*}
    Note that $\text{det}(P_{\epsilon}) = \frac{d_3}{\epsilon}-\Big( \frac{L+d_4}{2} \Big)^2$ and hence by taking $\epsilon < \epsilon_* := \frac{4d_3}{ (L+d_4)^2 }$, $P_{\epsilon}$ is negative definite. Therefore, by taking $\epsilon\in(0,\epsilon_*)$ and $\gamma\geq \gamma_*$, $\frac{d}{dt}V(\xi(t),\upsilon(t)) < 0$ whenever $(\xi(t),\upsilon(t))\in\Gamma$, unless $g_{\epsilon,\alpha,\gamma}(\xi,\upsilon)=0$ and $\xi-w(\upsilon)=0$. Now, for any point $(x,u)\in \Gamma$, such that $x=w(u)$, we have that $h(x) = \frac{1}{\beta^r}h_r(x,u)>0$. Thus, if $g_{\epsilon,\alpha,\gamma}(x,u)=0$ and $x=w(u)$,  from Prop. \ref{prop:critical_points_equilibria} we have that $(x,u)$ is a critical point of \eqref{eq:steady-state_opti_problem}; since $(x_*,u_*)$ is the sole critical point of \eqref{eq:steady-state_opti_problem}, $\frac{d}{dt}V(\xi(t),\upsilon(t)) < 0$ for all $(\xi(t),\upsilon(t))\in\Gamma\setminus\{(x_*,u_*)\}$, and $\frac{d}{dt}V(x_*,u_*) =0$. This completes the proof.
\end{proof}
\bibliography{mybib.bib}
\bibliographystyle{IEEEtran}

\end{document}